\documentclass[axioms,article,moreauthors,pdftex]{my-mdpi} 

\firstpage{1} 
\pubyear{2020}

\history{Submitted: September 9, 2020; Revised: September 30, October 14 and 16, 2020; Accepted: October 22, 2020}

\Title{Distributed-Order Non-Local Optimal Control}

\Author{Fa\"{\i}\c{c}al Nda\"{\i}rou $^{\dagger,\ddagger}$\orcidA{} 
and Delfim F. M. Torres $^{\ddagger}$*\orcidB{}}

\AuthorNames{Fa\"{\i}\c{c}al Nda\"{\i}rou and Delfim F. M. Torres}

\address[1]{Center for Research and Development in Mathematics and Applications (CIDMA),
Department of Mathematics, University of Aveiro, 3810-193 Aveiro, Portugal;
faical@ua.pt (F.N.); delfim@ua.pt (D.F.M.T.)}

\corres{Correspondence: delfim@ua.pt; Tel.: +351-234-370-668}

\firstnote{This research is part of first author's Ph.D. project, 
which is carried out at the University of Aveiro under 
the Doctoral Program in Applied Mathematics
of Universities of Minho, Aveiro, and Porto (MAP-PDMA).} 

\secondnote{These authors contributed equally to this work.}

\abstract{Distributed-order fractional non-local operators have been introduced 
and studied by Caputo at the end of the 20th century. They generalize 
fractional order derivatives/integrals in the sense that such operators 
are defined by a weighted integral of different orders of differentiation 
over a certain range. The subject of distributed-order non-local derivatives 
is currently under strong development due to its applications
in modeling some complex real world phenomena. Fractional optimal control theory 
deals with the optimization of a performance index functional subject 
to a fractional control system. One of the most important results in classical
and fractional optimal control is the Pontryagin Maximum Principle, 
which gives a necessary optimality condition that every solution 
to the optimization problem must verify. In our work, 
we extend the fractional optimal control theory by considering dynamical 
systems constraints depending on distributed-order fractional derivatives. 
Precisely, we prove a weak version of Pontryagin's maximum principle
and a sufficient optimality condition under appropriate convexity assumptions.}

\keyword{distributed-order fractional calculus;
basic optimal control problem; Pontryagin extremals.}

\MSC{26A33; 49K15}

\begin{document}
	

\section{Introduction}

Distributed-order fractional operators were introduced and studied by Caputo 
at the end of the previous century \cite{caputo1,caputo2}. They can be seen  
as a kind of generalization of fractional order derivatives/integrals 
in the sense that these operators are defined by a weighted integral 
of different orders of differentiation over a certain range. 
This subject gained more interest at the beginning of the current century, 
by researchers from different mathematical disciplines, through 
attempts to solve differential equations with distributed-order derivatives 
\cite{bagley1,bagley2,caputo3,diethelm}. Moreover, at the same time, 
in the domain of applied mathematics, those distributed-order fractional operators 
have started to be used, in a satisfactory way, to describe some complex phenomena 
modeling real world problem: see, for instance, works in viscoelasticity 
\cite{atan,lorenzo} and in diffusion \cite{antonio}. Today, the study
of distributed-order systems with fractional derivatives is a hot subject:
see, e.g., \cite{MR4151429,MR4109340,MR4046716} and references therein.

Fractional optimal control deals with optimization problems involving 
fractional differential equations as well as a performance index functional. 
One of the most important results is the Pontryagin Maximum Principle, 
which gives a first-order necessary optimality condition that every 
solution to the dynamic optimization problem must verify. By applying such result, 
it is possible to find and identify candidate solutions to the optimal control problem. 
For the state of the art on fractional optimal control we refer the readers
to \cite{MR3904404,MR3955185,MR3988048} and references therein.
Recently, distributed-order fractional problems of the calculus of variations
were introduced and investigated in \cite{ricardo}. Here, our main aim is to extend   
the distributed-order fractional Euler--Lagrange equation of \cite{ricardo}
to the Pontryagin setting (see Remark~\ref{rem:cor:rr}).

Regarding optimal control for problems with distributed-order fractional operators, 
the results are rare and reduce to the following two papers: \cite{MR3654793}
and \cite{zbMATH06915310}. Both works develop numerical methods while, in contrast, 
here we are interested in analytical results (not in numerical approaches). 
Moreover, our results are new and bring new insights. Indeed,
in \cite{MR3654793} the problem is considered with Riemann--Liouville 
distributed derivatives, while in our case we consider 
optimal control problems with Caputo distributed derivatives.
It should be also noted an inconsistence in \cite{MR3654793}:
when one defines the control system with a Riemann--Liouville derivative, 
then in the adjoint system it should appear a Caputo derivative; 
when one considers optimal control problems with  a control system with Caputo derivatives, 
then the adjoint equation should involve a Riemann--Liouville operator;
as a consequence of integration by parts (cf. Lemma~\ref{lemma1}).
This inconsistence has been corrected in \cite{zbMATH06915310}, where
optimal control problems with Caputo distributed derivatives (like we do here) are considered.
Unfortunately, there is still an inconsistence in the necessary optimality conditions 
of both \cite{MR3654793} and \cite{zbMATH06915310}: the transversality conditions are written there
exactly as in the classical case, with the multiplier vanishing at the end of the interval, 
while the correct condition, as we prove in our Theorem~\ref{theo}, should involve
a distributed integral operator: see condition \eqref{trans}.

The text is organized as follows.
We begin by recalling definitions and necessary
results of the literature in Section~\ref{sec:2:P} of preliminaries. 
Our original results are then given in Section~\ref{sectionR}.
More precisely, we consider fractional optimal control problems where 
the dynamical system constraints depend on distributed-order fractional derivatives. 
We prove a weak version of Pontryagin's maximum principle 
for the considered distributed-order fractional problems (see Theorem~\ref{theo}) 
and investigate a Mangasarian-type sufficient optimality condition
(see Theorem~\ref{theosuff}). An example, illustrating the usefulness
of the obtained results, is given (see Examples~\ref{propos} and \ref{propost}).
We end with Section~\ref{sec:conc} of conclusions, mentioning also some possibilities
of future research.


\section{Preliminaries}
\label{sec:2:P}

In this section, we recall necessary results and fix notations.
We assume the reader to be familiar with the standard
Riemann--Liouville and Caputo fractional calculi 
\cite{MR3443073,MR1347689}. 

Let $\alpha$ be a real number in $[0, 1]$ and let $\psi$ 
be a non-negative continuous function defined on $[0, 1]$ 
such that 
$$
\int^1_0 \psi(\alpha)d\alpha >0.
$$ 
This function $\psi$ will act as a distribution of the order of differentiation.

\begin{Definition}[See \cite{caputo1}]
The left and right-sided Riemann--Liouville distributed-order 
fractional derivatives of a function 
$x: [a, b]\rightarrow \mathbb{R}$ are defined, respectively, by 
\[
\mathbb{D}^{\psi(\cdot)}_{a^{+}}x(t)= \int^1_0 \psi(\alpha)
\cdot D^{\alpha}_{a^{+}}x(t)d\alpha 
\quad \text{ and } \quad 
\mathbb{D}^{\psi(\cdot)}_{b^{-}}x(t)
= \int^1_0 \psi(\alpha)\cdot D^{\alpha}_{b^{-}}x(t)d\alpha,
\]
where $D^{\alpha}_{a^{+}}$ and $D^{\alpha}_{b^{-}}$ are, 
respectively, the left and right-sided Riemann--Liouville 
fractional derivatives of order $\alpha$.
\end{Definition}

\begin{Definition}[See \cite{caputo1}]
The left and right-sided Caputo distributed-order fractional 
derivatives of a function $x: [a, b]\rightarrow \mathbb{R}$ 
are defined, respectively, by
\[
^{C}\mathbb{D}^{\psi(\cdot)}_{a^{+}}x(t)= \int^1_0 \psi(\alpha)
\cdot ^{C}D^{\alpha}_{a^{+}}x(t)d\alpha 
\quad \text{ and } \quad  
^{C}\mathbb{D}^{\psi(\cdot)}_{b^{-}}x(t)
= \int^1_0 \psi(\alpha)\cdot ^{C}D^{\alpha}_{b^{-}}x(t)d\alpha,
\]
where $^{C}D^{\alpha}_{a^{+}}$ and $^{C}D^{\alpha}_{b^{-}}$ are,
respectively, the left and right-sided Caputo fractional 
derivatives of order $\alpha$.
\end{Definition}

As noted in \cite{ricardo}, there is a relation between the 
Riemann--Liouville and the Caputo distributed-order fractional derivatives:
\[
^{C}\mathbb{D}^{\psi(\cdot)}_{a^{+}}x(t)
= \mathbb{D}^{\psi(\cdot)}_{a^{+}}x(t)
- x(a)\int^1_0 \frac{\psi(\alpha)}{\Gamma(1-\alpha)}(t-a)^{-\alpha}d\alpha
\]
and 
\[
^{C}\mathbb{D}^{\psi(\cdot)}_{b^{-}}x(t)
= \mathbb{D}^{\psi(\cdot)}_{b^{-}}x(t)- x(b)
\int^1_0 \frac{\psi(\alpha)}{\Gamma(1-\alpha)}(b-t)^{-\alpha}d\alpha.
\]

Along the text, we use the notation
\[
\mathbb{I}^{1-\psi(\cdot)}_{b^{-}}x(t)
=\int^1_0\psi(\alpha)\cdot I^{1-\alpha}_{b^{-}}x(t)d\alpha,
\]
where $I^{1-\alpha}_{b^{-}}$ represents the right 
Riemann--Liouville fractional integral of order $1-\alpha$.

The next result has an essential role in the proofs
of our main results, that is, in the proofs of 
Theorems~\ref{theo} and \ref{theosuff}.

\begin{Lemma}[Integration by parts formula \cite{ricardo}]
\label{lemma1}
Let $x$ be a continuous function and $y$ a continuously 
differentiable function. Then, 
\[
\int^b_a x(t)\cdot ^{C}\mathbb{D}^{\psi(\cdot)}_{a^{+}}y(t)dt 
= \left[ y(t)\cdot \mathbb{I}^{1-\psi (\cdot)}_{b^{-}}x(t) \right]^b_a 
+ \int^b_a y(t)\cdot \mathbb{D}^{\psi(\cdot)}_{b^{-}}x (t)dt.
\]
\end{Lemma}

Next, we recall the standard notion of concave function, 
which will be used in Section~\ref{subsec:SCGO}.

\begin{Definition}[See \cite{MR0274683}]
A function $h: \mathbb{R}^{n} \rightarrow \mathbb{R}$ is concave if 
\[
h(\beta \theta_1 + (1- \beta)\theta_2)
\geq \beta h(\theta_1)+ (1-\beta)h(\theta_2)
\]
for all $\beta \in [0, 1]$ and for all 
$\theta_1$, $\theta_2 $ in $\mathbb{R}^n$. 
\end{Definition}

\begin{Lemma}[See \cite{MR0274683}]
\label{lemma:concave}
Let $h: \mathbb{R}^{n} \rightarrow \mathbb{R}$ be a continuously differentiable function. 
Then $h$ is a concave function if and only if it satisfies the so called gradient inequality:
\[
h(\theta_1)-h(\theta_2)\geq \nabla h(\theta_1)(\theta_1- \theta_2)
\]
for all $\theta_1, \theta_2 \in \mathbb{R}^n$.
\end{Lemma}

Finally, we recall a fractional version of Gronwall's inequality,
which will be useful to prove continuity of solutions 
in Section~\ref{sub:sec:CS}.

\begin{Lemma}[See \cite{gronwall}]
\label{thm:gronwall}
Let $\alpha$ be a positive real number and let $a(\cdot)$, 
$b(\cdot)$, and $u(\cdot)$ be non-negative continuous functions 
on $[0, T]$ with $b(\cdot)$  monotonic increasing on $[0, T)$.
If
\[
u(t)\leq a(t) + b(t)\int^t_0(t-s)^{\alpha-1}u(s)ds,
\]
then 
\[
u(t)\leq a(t) + \int^t_0 \left[ \sum^{\infty}_{n=0}
\frac{\left(b(t)\Gamma(\alpha ) \right)^n}{\Gamma(n\alpha)} 
(t-s)^{n\alpha-1}u(s)\right]ds
\]
for all $t\in [0,T)$.
\end{Lemma}


\section{Main Results}
\label{sectionR}

The basic problem of optimal control we consider in this work, 
denoted by \eqref{bp}, consists to find a piecewise continuous 
control $u \in PC$ and the corresponding piecewise smooth 
state trajectory $x \in PC^{1}$, solution of the  
distributed-order non-local variational problem
\begin{equation}
\tag{BP}
\label{bp}
\begin{gathered}
J[x(\cdot), u(\cdot)]= \int^{b}_{a} L\left(t, x(t), u(t)\right)dt 
\longrightarrow \max,\\
^{C}\mathbb{D}^{\psi(\cdot)}_{{a}^{+}}x(t)
= f\left(t, x(t), u(t)\right), \quad t\in [a, b],\\
x(\cdot) \in PC^{1}, \quad u(\cdot) \in PC,\\
x(a)= x_a,
\end{gathered}
\end{equation}
where functions $L$ and $f$, both defined on $[a, b]\times 
\mathbb{R}\times \mathbb{R}$, are assumed to be continuously differentiable 
in all their three arguments: $L \in C^1$, $f\in C^1$. Our main
contribution is to prove necessary (Section~\ref{subsec:PMP})
and sufficient (Section~\ref{subsec:SCGO}) optimality conditions.


\subsection{Sensitivity analysis} 
\label{sub:sec:CS}

Before we can prove necessary optimality conditions to problem \eqref{bp},
we need to establish continuity and differentiability results on the state solutions 
for any control perturbation (Lemmas~\ref{lemma2} and \ref{lemma3}), which are then
used in Section~\ref{subsec:PMP}. The proof of Lemma~\ref{lemma2} makes
use of the following mean value theorem for integration, that can be found
in any textbook of calculus (see, e.g., Lemma~1 of \cite{MR3089388}):
if $F : [0,1] \rightarrow \mathbb{R}$ is a continuous function and $\psi$
is an integrable function that does not change sign on the interval, 
then there exists a number $\bar{\alpha}$ such that
$$
\int_0^1 \psi(\alpha) F(\alpha) d\alpha 
= F(\bar{\alpha}) \int_0^1 \psi(\alpha) d\alpha.
$$

\begin{Lemma}[Continuity of solutions]
\label{lemma2} 
Let $u^{\epsilon}$ be a control perturbation around the optimal control $u^{*}$, 
that is, for all $t\in [a,b]$, $ u^{\epsilon}(t)= u^{*}(t)+\epsilon h(t)$, 
where $h(\cdot) \in PC$ is a variation and $\epsilon \in \mathbb{R}$. 
Denote by $x^{\epsilon}$ its corresponding state trajectory, solution of 
\[
^{C}\mathbb{D}^{\psi(\cdot)}_{{a}^{+}} x^{\epsilon}(t)
= f\left(t, x^{\epsilon}(t), u^{\epsilon}(t)\right), \quad x^{\epsilon}(a)= x_a.
\]
Then, we have that $x^{\epsilon}$ converges to the optimal 
state trajectory $x^{*}$ when $\epsilon$ tends to zero.
\end{Lemma}

\begin{proof}
Starting from definition, we have, for all $t\in [a,b]$, that
\[
\left|^{C}\mathbb{D}^{\psi(\cdot)}_{{a}^{+}} x^{\epsilon}(t) 
- ^{C}\mathbb{D}^{\psi(\cdot)}_{{a}^{+}} x^{*}(t) \right|
= \left| f\left(t, x^{\epsilon}(t), u^{\epsilon}(t)\right) 
- f\left(t, x^{*}(t), u^{*}(t)\right)\right|.
\]
Then, by linearity,
\[
\left|^{C}\mathbb{D}^{\psi(\cdot)}_{{a}^{+}} x^{\epsilon}(t) 
- ^{C}\mathbb{D}^{\psi(\cdot)}_{{a}^{+}} x^{*}(t) \right|
= \left|^{C}\mathbb{D}^{\psi(\cdot)}_{{a}^{+}}\left( 
x^{\epsilon}(t)- x^{*}(t)\right) \right|
= \left| f\left(t, x^{\epsilon}(t), u^{\epsilon}(t)\right) 
- f\left(t, x^{*}(t), u^{*}(t)\right)\right|
\]
and it follows, by definition of the distributed operator, that
\begin{equation*}
\left|\int^1_0 \psi(\alpha)^CD^{\alpha}_{{a}^{+}}\left( x^{\epsilon}(t)
- x^{*}(t)\right)  d\alpha \right|
= \left| f\left(t, x^{\epsilon}(t), u^{\epsilon}(t)\right) 
- f\left(t, x^{*}(t), u^{*}(t)\right)\right|.
\end{equation*}
Now, using the mean value theorem for integration,
and denoting $m:= \int_{0}^{1} \psi(\alpha) d\alpha$, 
we obtain that there exists an $\bar{\alpha}$ such that
\begin{equation*}
\left|
{^CD^{\bar{\alpha}}_{{a}^{+}}}\left( x^{\epsilon}(t)
- x^{*}(t)\right) \right|
\leq \frac{\left| f\left(t, x^{\epsilon}(t), u^{\epsilon}(t)\right) 
- f\left(t, x^{*}(t), u^{*}(t)\right)\right|}{m}.
\end{equation*}
Clearly, one has
\begin{equation*}
^CD^{\bar{\alpha}}_{{a}^{+}}\left( x^{\epsilon}(t)
- x^{*}(t)\right) 
\leq \left|^CD^{\bar{\alpha}}_{{a}^{+}}\left( x^{\epsilon}(t)
- x^{*}(t)\right) \right| 
\leq \frac{\left| f\left(t, x^{\epsilon}(t), u^{\epsilon}(t)\right) 
- f\left(t, x^{*}(t), u^{*}(t)\right)\right|}{m},
\end{equation*}
which leads to  
\begin{equation*}
x^{\epsilon}(t)- x^{*}(t)
\leq
I^{\bar{\alpha}}_{a^{+}} \left[
\frac{\left| f\left(t, x^{\epsilon}(t), u^{\epsilon}(t)\right) 
- f\left(t, x^{*}(t), u^{*}(t)\right)\right|}{m}
\right].
\end{equation*}
Moreover, because $f$ is Lipschitz-continuous, we have
\[
\Bigl| f\left(t, x^{\epsilon}, u^{\epsilon}\right) 
- f\left(t, x^{*}, u^{*}\right)\Bigr| \leq K_1 \Bigl| 
x^{\epsilon}- x^{*}\Bigr| + K_2\bigl|u^{\epsilon}-u^{*}
\Bigr|. 
\]
By setting $K= \max\{K_1, K_2\}$, it follows that
\begin{align*}
\Bigl| x^{\epsilon}(t) - x^{*}(t) \Bigr|
&\leq \frac{K}{m}I^{\bar{\alpha}}_{a^{+}}\Big( \Bigl| x^{\epsilon}(t)
- x^{*}(t)\Bigr| + \Bigl|\epsilon h(t)\Bigr| \Big)\\
&= \frac{K}{m}\left[ |\epsilon|I^{\bar{\alpha}}_{a^{+}} 
\left(\Bigl| h(t) \Bigr| \right)  + I^{\bar{\alpha}}_{a^{+}} 
\left( \Bigl|x^{\epsilon}(t)- x^{*}(t) \Bigr| \right)  \right]\\
&= \frac{K}{m}\left[ |\epsilon|I^{\bar{\alpha}}_{a^{+}} \left(\Bigl| h(t) \Bigr| \right)  
+ \frac{1}{\Gamma(\bar{\alpha})}\int^t_a(t-s)^{\bar{\alpha} -1}
\Bigl|x^{\epsilon}(s)- x^{*}(s) \Bigr|ds \right]
\end{align*}
for all $t\in [a,b]$. Now, by applying Lemma~\ref{thm:gronwall} 
(the fractional Gronwall inequality), it follows that
\begin{align*}
\Bigl| x^{\epsilon}(t)- x^{*}(t) \Bigr|
&\leq \frac{K}{m}\left[ |\epsilon|I^{\bar{\alpha}}_{a^{+}} 
\left( \Bigl| h(t) \Bigr| \right)  + |\epsilon|
\int^t_a \left( \sum^{\infty}_{i=0} \frac{1}{\Gamma(i\bar{\alpha})}(t-s)^{i\bar{\alpha} -1} 
I^{\bar{\alpha}}_{a^{+}} \left( \Bigl| h(s) \Bigr| \right) \right) ds \right]\\
&= |\epsilon|\frac{K}{m}\left[ I^{\bar{\alpha}}_{a^{+}} \left( \Bigl| h(t) \Bigr| \right)  
+ \int^t_a \left( \sum^{\infty}_{i=1} \frac{1}{\Gamma(i\bar{\alpha} + 1)}(t-s)^{i\bar{\alpha}} 
I^{\bar{\alpha}}_{a^{+}} \left( \Bigl| h(s) \Bigr| \right) \right) ds \right]\\
&\leq |\epsilon|\frac{K}{m}\left[ I^{\bar{\alpha}}_{a^{+}} \left( \Bigl| h(t) \Bigr| \right)  
+ \int^t_a \left( \sum^{\infty}_{i=1} \frac{\delta^{i\bar{\alpha}}}{\Gamma(i\bar{\alpha} + 1)} 
I^{\bar{\alpha}}_{a^{+}} \left( \Bigl| h(s) \Bigr| \right) \right) ds \right].
\end{align*}
The series in the last inequality is a Mittag--Leffler function and thus convergent. 
Hence, by taking the limit when $\epsilon$ tends to zero, we obtain the desired result: 
$x^{\epsilon} \rightarrow x^{*}$ for all $t\in [a, b]$.
\end{proof}

\begin{Lemma}[Differentiability of the perturbed trajectory]
\label{lemma3}
There exists a function $\eta$ defined on $[a,b]$ such that 
\[
x^{\epsilon}(t)= x^{*}(t) + \epsilon \eta(t) + o(\epsilon).
\]
\end{Lemma} 

\begin{proof}
Since $f \in C^1$, we have that
\begin{equation*}
f(t, x^{\epsilon}, u^{\epsilon})
= f(t, x^{*}, u^{*}) + (x^{\epsilon}-x^{*})
\frac{\partial f(t, x^{*}, u^{*})}{\partial x} 
+ (u^{\epsilon}-u^{*})\frac{\partial f(t, x^{*}, u^{*})}{\partial u}
+ o(|x^{\epsilon} - x^{*}|,|u^{\epsilon}-u^{*}|).
\end{equation*}
Observe that $u^{\epsilon}-u^{*}= \epsilon h(t)$ 
and $u^{\epsilon} \rightarrow u^{*}$ when $\epsilon \rightarrow 0$ and, 
by Lemma~\ref{lemma2}, we have $x^{\epsilon} \rightarrow x^{*}$ 
when $\epsilon \rightarrow 0$. Thus, the residue term can be expressed 
in terms of $\epsilon$ only, that is, the residue is $o(\epsilon)$. 
Therefore, we have
\[
^C\mathbb{D}^{\psi(\cdot)}_{a^{+}} x^{\epsilon}(t)
= ^C\mathbb{D}^{\psi(\cdot)}_{a^{+}} x^{*}(t) 
+ (x^{\epsilon}-x^{*})\frac{\partial f(t, x^{*}, u^{*})}{\partial x} 
+ \epsilon h(t)\frac{\partial f(t, x^{*}, u^{*})}{\partial u} + o(\epsilon),
\]
which leads to
\[
\lim_{\epsilon \rightarrow 0} \left[\frac{^C\mathbb{D}^{\psi(\cdot)}_{a^{+}}
(x^{\epsilon}-x^{*})}{\epsilon}-\frac{(x^{\epsilon}-x^{*})}{\epsilon}
\frac{\partial f(t, x^{*}, u^{*})}{\partial x} - h(t)
\frac{\partial f(t, x^{*}, u^{*})}{\partial u}\right]=0,
\]
meaning that
\[
^C\mathbb{D}^{\psi(\cdot)}_{a^{+}}\Big(\lim_{\epsilon \rightarrow 0} 
\frac{x^{\epsilon}-x^{*}}{\epsilon} \Big) = \lim_{\epsilon \rightarrow 0}
\frac{x^{\epsilon}-x^{*}}{\epsilon}\frac{\partial 
f(t, x^{*}, u^{*})}{\partial x} + h(t)\frac{\partial 
f(t, x^{*}, u^{*})}{\partial u}.
\]
We want to prove the existence of the limit
$\displaystyle{\lim_{\epsilon \rightarrow 0} 
\frac{x^{\epsilon}-x^{*}}{\epsilon}}=:\eta$, that is, 
to prove that 
$x^{\epsilon}(t)= x^{*}(t) + \epsilon \eta(t) + o(\epsilon)$.
This is indeed the case, since $\eta$ is solution of the distributed 
order fractional differential equation
\begin{equation*}
\begin{cases}
^C\mathbb{D}^{\psi(\cdot)}_{a^{+}} \eta(t)
= \frac{\partial f(t, x^{*}, u^{*})}{\partial x}\eta(t) 
+\frac{\partial f(t, x^{*}, u^{*})}{\partial u}h(t),\\[3mm]
\eta(a)= 0.
\end{cases}
\end{equation*}
The intended result is proved.
\end{proof}


\subsection{Pontryagin's maximum principle of distributed-order}
\label{subsec:PMP}

The following result is a necessary condition 
of Pontryagin type \cite{MR0186436} for the basic 
distributed-order non-local optimal control problem \eqref{bp}.

\begin{Theorem}[Pontryagin Maximum Principle for \eqref{bp}]
\label{theo}
If $(x^{*}(\cdot), u^{*}(\cdot))$ is an optimal pair for \eqref{bp}, 
then there exists $\lambda \in PC^1$, called the adjoint function variable, 
such that the following conditions hold for all $t$ in the interval $[a, b]$:
\begin{itemize}

\item the optimality condition
\begin{equation}
\label{opt}
\frac{\partial L}{\partial u}(t, x^{*}(t), u^{*}(t))
+ \lambda (t)\frac{\partial f}{\partial u}(t, x^{*}(t), u^{*}(t))=0;
\end{equation}
		
\item the adjoint equation
\begin{equation}
\label{adj}
\mathbb{D}^{\psi(\cdot)}_{b^{-}} \lambda(t)
= \frac{\partial L}{\partial x}(t,x^{*}(t), u^{*}(t)) 
+ \lambda (t)\frac{\partial f}{\partial x }(t, x^{*}(t), u^{*}(t));
\end{equation}
		
\item the transversality condition
\begin{equation}
\label{trans}
\mathbb{I}^{1-\psi(\cdot)}_{b^{-}}\lambda(b)=0.
\end{equation}
\end{itemize}
\end{Theorem}

\begin{proof}
Let $(x^{*}(\cdot), u^{*}(\cdot))$ be solution to problem \eqref{bp}, 
$h(\cdot) \in PC$ be a variation, and $\epsilon$ a real constant. 
Define $u^{\epsilon}(t)= u^{*}(t)+ \epsilon h(t)$, so that 
$u^{\epsilon}\in PC$. Let $x^{\epsilon}$ be the state corresponding 
to the control $u^{*}$, that is, the state solution of 
\begin{equation}
\label{equatepsi}
^{C}\mathbb{D}^{\psi(\cdot)}_{{a}^{+}} x^{\epsilon}(t)
= f\left(t, x^{\epsilon}(t), u^{\epsilon}(t)\right), 
\quad x^{\epsilon}(a)= x_a.
\end{equation}
Note that $u^{\epsilon}(t) \rightarrow u^{*}(t)$ for all 
$t\in [a, b]$ whenever $\epsilon \rightarrow 0$. Furthermore, 
\begin{equation}
\label{eqpartialu}
\frac{\partial u^{\epsilon}(t)}{\partial \epsilon}\Bigr|_{\epsilon=0}  = h(t).
\end{equation} 
Something similar is also true for $x^{\epsilon}$: because $ f\in C^1$,
it follows from Lemma~\ref{lemma2} that, for each fixed $t$, 
$x^{\epsilon}(t)\rightarrow x^{*}(t)$ as $\epsilon \rightarrow 0$. 
Moreover, by Lemma~\ref{lemma3}, the derivative 
$\displaystyle{\frac{\partial x^{\epsilon}(t)}{\partial \epsilon}\Bigr|_{\epsilon=0}}$ 
exists for each $t$. The objective functional at $(x^{\epsilon}, u^{\epsilon})$ is 
\[
J[x^{\epsilon}, u^{\epsilon}]
= \int^b_a L\left(t, x^{\epsilon}(t), u^{\epsilon}(t) \right)dt.
\]
Next, we introduce the adjoint function $\lambda$. 
Let $\lambda(\cdot)$ be in $PC^1$, to be determined. 
By the integration by parts formula (see Lemma~\ref{lemma1}),
\[
\int^b_a \lambda(t)\cdot ^{C}\mathbb{D}^{\psi(\cdot)}_{a^{+}}
x^{\epsilon}(t)dt = \left[ x^{\epsilon}(t)\cdot \mathbb{I}^{1
-\psi (\cdot)}_{b^{-}}\lambda(t) \right]^b_a + \int^b_a
x^{\epsilon}(t)\cdot \mathbb{D}^{\psi(\cdot)}_{b^{-}}\lambda (t)dt,
\]
and one has
\begin{equation*}
\int^b_a \lambda(t)\cdot ^{C}\mathbb{D}^{\psi(\cdot)}_{a^{+}}x^{\epsilon}(t)dt 
- \int^b_ax^{\epsilon}(t)\cdot \mathbb{D}^{\psi(\cdot)}_{b^{-}}\lambda (t)dt 
- x^{\epsilon}(b)\cdot \mathbb{I}^{1-\psi (\cdot)}_{b^{-}}\lambda(b)
+ x^{\epsilon}(a)\cdot \mathbb{I}^{1-\psi (\cdot)}_{b^{-}}\lambda(a) =0.
\end{equation*}
Adding this zero to the expression $J[x^{\epsilon}, u^{\epsilon}]$ gives 
\begin{multline*}
\phi (\epsilon)= J[x^{\epsilon}, u^{\epsilon}]
= \int^b_a \left[L\left( t, x^{\epsilon}(t), u^{\epsilon}(t)\right) 
+ \lambda(t)\cdot ^{C}\mathbb{D}^{\psi(\cdot)}_{a^{+}}x^{\epsilon}(t)
- x^{\epsilon}(t)\cdot \mathbb{D}^{\psi(\cdot)}_{b^{-}}\lambda (t) \right]dt\\
- x^{\epsilon}(b)\cdot \mathbb{I}^{1-\psi (\cdot)}_{b^{-}}\lambda(b) 
+ x^{\epsilon}(a)\cdot \mathbb{I}^{1-\psi (\cdot)}_{b^{-}}\lambda(a), 
\end{multline*}
which by \eqref{equatepsi} is equivalent to
\begin{multline*}
\phi (\epsilon)= J[x^{\epsilon}, u^{\epsilon}]
= \int^b_a \left[ L\left( t, x^{\epsilon}(t), u^{\epsilon}(t)\right) 
+ \lambda(t)\cdot f\left(t, x^{\epsilon}(t), u^{\epsilon}(t) \right)
- x^{\epsilon}(t)\cdot \mathbb{D}^{\psi(\cdot)}_{b^{-}}
\lambda (t) \right]dt\\ 
-  x^{\epsilon}(b)\cdot \mathbb{I}^{1
-\psi(\cdot)}_{b^{-}}\lambda(b) + x_a\cdot \mathbb{I}^{1
-\psi (\cdot)}_{b^{-}}\lambda(a). 
\end{multline*}
Since the process $(x^{*}, u^{*})= (x^0, u^0)$ is assumed to be 
a maximizer of problem \eqref{bp}, the derivative of 
$\phi(\epsilon)$ with respect to $\epsilon$  
must vanish at $\epsilon=0 $, that is,
\begin{equation*}
\begin{split}
0&= \phi'(0)
= \frac{d }{d \epsilon } J[x^{\epsilon}, u^{\epsilon}]|_{\epsilon=0}\\
&= \int^b_a \left[ \frac{\partial L}{\partial x}
\frac{\partial x^{\epsilon}(t)}{\partial \epsilon}\Bigr|_{\epsilon=0} 
+ \frac{\partial L}{\partial u}\frac{\partial u^{\epsilon}(t)}{\partial 
\epsilon}\Bigr|_{\epsilon=0} 
+ \lambda (t)\left( \frac{\partial f}{\partial x}
\frac{\partial x^{\epsilon}(t)}{\partial \epsilon}\Bigr|_{\epsilon=0}  
+ \frac{\partial f}{\partial u}\frac{\partial u^{\epsilon}(t)}{\partial 
\epsilon}\Bigr|_{\epsilon=0}\right) \right. \\
&\qquad\qquad  \left. - \mathbb{D}^{\psi(\cdot)}_{b^{-}}\lambda(t)
\frac{\partial x^{\epsilon}(t)}{\partial \epsilon}\Bigr|_{\epsilon=0}\right]dt
-\frac{\partial x^{\epsilon}(b)}{\partial \epsilon}\Bigr|_{\epsilon=0} 
\mathbb{I}^{1-\psi (\cdot)}_{b^{-}}\lambda(b),
\end{split}
\end{equation*}
where the partial derivatives of $L$ and $f$, with respect to $x$ and $u$, 
are evaluated at $\left( t, x^{*}(t), u^{*}(t) \right)$. Rearranging 
the term and using \eqref{eqpartialu}, we obtain that
\begin{equation*}
\int^b_a \left[ \Big(\frac{\partial L}{\partial x} 
+ \lambda(t)\frac{\partial f}{\partial x } 
- \mathbb{D}^{\psi(\cdot)}_{b^{-}}\lambda(t)\Big) 
\frac{\partial x^{\epsilon}(t)}{\partial \epsilon}\Bigr|_{\epsilon=0} 
+ \Big(\frac{\partial L}{\partial u} + \lambda(t)\frac{\partial 
f}{\partial u}\Big)h(t)\right]dt
-\frac{\partial x^{\epsilon}(b)}{\partial 
\epsilon}\Bigr|_{\epsilon=0} \mathbb{I}^{1-\psi (\cdot)}_{b^{-}}\lambda(b)=0.
\end{equation*}
Setting $H(t,x,u,\lambda)= L(t,x,u) + \lambda f(t,x,u)$, it follows that
\begin{equation*}
\int^b_a \left[ \Big(\frac{\partial H}{\partial x}  
- \mathbb{D}^{\psi(\cdot)}_{b^{-}}\lambda(t)\Big) 
\frac{\partial x^{\epsilon}(t)}{\partial \epsilon}\Bigr|_{\epsilon=0} 
+ \frac{\partial H}{\partial u} h(t)\right]dt \\
-\frac{\partial x^{\epsilon}(b)}{\partial \epsilon}\Bigr|_{\epsilon=0} 
\mathbb{I}^{1-\psi (\cdot)}_{b^{-}}\lambda(b)=0,
\end{equation*}
where the partial derivatives of $H$ are evaluated at
$\left( t, x^{*}(t), u^{*}(t), \lambda(t) \right)$. Now, choosing
\[
\mathbb{D}^{\psi(\cdot)}_{b^{-}}\lambda(t) 
= \frac{\partial H}{\partial x}\left( t, x^{*}(t), u^{*}(t), \lambda(t) \right), 
\quad \text{ with } \mathbb{I}^{1-\psi (\cdot)}_{b^{-}}
\lambda(b)=0,
\]
that is, given the adjoint equation \eqref{adj} 
and the transversality condition \eqref{trans}, 
it yields 
\[
\int^b_a \frac{\partial H}{\partial 
u}\left(t, x^{*}(t), u^{*}(t), \lambda(t) \right) h(t)=0
\]
and, by the fundamental lemma of the calculus of variations 
\cite{MR500859}, we have the optimality condition \eqref{opt}:
\[
\frac{\partial H}{\partial u}\left(t, x^{*}(t), u^{*}(t), \lambda(t) \right)=0.
\]
This concludes the proof.
\end{proof}

\begin{Remark}
If we change the basic optimal control problem \eqref{bp} by changing
the boundary condition given on the state variable
at initial time, $x(a)= x_a$, to a terminal condition, then 
the optimality condition and the adjoint equation 
of the Pontryagin Maximum Principle (Theorem~\ref{theo}) 
remain exactly the same. Changes appear only
on the transversality condition:
\begin{itemize}
\item a boundary condition at final/terminal time, that is,
fixing the value $x(b)=x_b$ with $x(a)$ remaining free, 
leads to
$$ 
\mathbb{I}^{1-\psi(\cdot)}_{a^{-}}\lambda(a)=0; 
$$
\item in case when no boundary conditions is given
(i.e., both $x(a)$ and $x(b)$ are free), then we have
$$ 
\mathbb{I}^{1-\psi(\cdot)}_{b^{-}}\lambda(b)=0 
\quad \text{ and } \quad  
\mathbb{I}^{1-\psi(\cdot)}_{a^{-}}\lambda(a)=0.
$$ 
\end{itemize}
\end{Remark}

\begin{Remark}
\label{rem:cor:rr}
If $f\left(t, x, u\right) = u$, that is, $^{C}{\mathbb{D}}^{\psi(\cdot)}_{a+} x(t)= u(t)$, 
then our problem \eqref{bp} gives a basic problem of the calculus of variations, 
in the distributed-order fractional sense of \cite{ricardo}. In this very particular case, 
we obtain from our Theorem~\ref{theo} the Euler--Lagrange equation of \cite{ricardo} 
(cf. Theorem~2 of \cite{ricardo}).
\end{Remark}

\begin{Remark}
Our distributed-order fractional optimal control problem \eqref{bp} 
can be easily extended to the vector setting. Precisely, let 
$x:= (x_1, \ldots , x_n)$ and $u:= (u_1, \ldots , u_n)$ 
with $(n, m) \in \mathbb{N}^2$ such that $m\leq n$, 
and functions $f: [a, b]\times \mathbb{R}^n\times \mathbb{R}^m 
\rightarrow \mathbb{R}^n$ and 
$ L: [a, b]\times \mathbb{R}^n\times \mathbb{R}^m \rightarrow \mathbb{R}$ 
be  continuously differentiable with respect to all its components. 
If $(x^{*}, u^{*})$ is an optimal pair, then the following 
conditions hold for $t\in [a, b]$:
\begin{itemize}
\item the optimality conditions
\begin{equation*}
\frac{\partial L}{\partial u_i}(t, x^{*}(t), u^{*}(t))
+ \lambda(t) \cdot \frac{\partial f}{\partial u_i} (t, x^{*}(t), u^{*}(t))=0, 
\quad i=1, \ldots , m;
\end{equation*}
\item the adjoint equations
\begin{equation*}
\mathbb{D}^{\psi(\cdot)}_{b^{-}} \lambda_j(t)
= \frac{\partial L}{\partial x_j}(t,x^{*}(t), u^{*}(t)) 
+ \lambda(t) \cdot \frac{\partial f}{\partial x_j }(t, x^{*}(t), u^{*}(t)), 
\quad j=1, \ldots, n;
\end{equation*}
\item the transversality conditions
\begin{equation}\label{trans:b}
\mathbb{I}^{1-\psi(\cdot)}_{b^{-}}\lambda_j(b)=0,
\quad j=1, \ldots, n.
\end{equation}
\end{itemize}
\end{Remark}

\begin{Definition}
The candidates to solutions of \eqref{bp}, 
obtained by the application of our Theorem~\ref{theo}, 
will be called (Pontryagin) extremals.
\end{Definition}

We now illustrate the usefulness of our Theorem~\ref{theo}
with an example.

\begin{Example}
\label{propos}
The triple $(\tilde{x}, \tilde{u}, \lambda)$ given by 
$\tilde{x}(t)= t^2$, $\tilde{u}(t)= \displaystyle{\frac{t(t-1)}{\ln t }}$, 
and $\lambda (t)=0$, for  $t\in [0,1]$, is an extremal 
of the following distributed-order fractional optimal control problem:
\begin{equation}
\label{pexmple}
\begin{gathered}
J[x(\cdot), u(\cdot)]= \int^{1}_{0} -\left(x(t)-t^2\right)^2 
- \left(u- \frac{t(t-1)}{\ln t} \right)^2 \longrightarrow \max,\\
^{C}\mathbb{D}^{\psi(\cdot)}_{{0}^{+}}x(t)= u(t), \quad t\in [0, 1],\\
x(0)= 0.
\end{gathered}
\end{equation}
Indeed, by defining the Hamiltonian function as 
\begin{equation}
\label{hami}
H(t, x, u, \lambda)= -\left[(x-t^2) 
+ \left(u- \frac{t(t-1)}{\ln t} \right)^2 \right] + \lambda u,
\end{equation}  
it follows:
\begin{itemize}
\item from  the optimality condition 
$\displaystyle{\frac{\partial H}{\partial u}=0}$,
\begin{equation}
\label{optiexmple}
\lambda(t)= 2 \left(u- \frac{t(t-1)}{\ln t} \right);
\end{equation} 
\item from the adjoint equation 
$\mathbb{D}^{\psi(\alpha)}_{{0}^{+}}\lambda(t)
= \displaystyle{\frac{\partial H}{\partial x}}$, 
\begin{equation}
\label{adjoexmpl}
\mathbb{D}^{\psi(\alpha)}_{{0}^{+}}\lambda(t)= -2(x-t^2);
\end{equation}
\item from the transversality condition, 
\begin{equation}
\label{transiexmpl}
\mathbb{I}^{1-\psi(\alpha)}_{b^{-}}\lambda(b)=0.
\end{equation}
\end{itemize}
We easily see that \eqref{optiexmple}, \eqref{adjoexmpl} 
and \eqref{transiexmpl} are satisfied for 
$$
x(t)= t^2, 
\quad u(t)= \frac{t(t-1)}{\ln t},
\quad  \lambda(t)=0.
$$
\end{Example}


\subsection{Sufficient condition for global optimality}
\label{subsec:SCGO}

We now prove a Mangasarian type theorem for the 
distributed-order fractional optimal control problem \eqref{bp}.

\begin{Theorem}
\label{theosuff}
Consider the basic distributed-order fractional optimal control problem \eqref{bp}. 
If $(x, u) \rightarrow L(t, x, u)$ and $(x, u) \rightarrow f(t, x, u)$ 
are concave and $(\tilde{x}, \tilde{u}, \lambda)$ is a Pontryagin extremal 
with $\lambda (t)\geq 0$, $t \in [a,b]$, then 
$$
J[\tilde{x}, \tilde{u}] \geq J[x,u] 
$$ 
for any admissible pair $(x,u)$.
\end{Theorem}

\begin{proof}
Because $L$ is concave as a function of $x$ and $u$, 
we have from Lemma~\ref{lemma:concave} that
\begin{equation*}
L\left(t, \tilde{x}(t), \tilde{u}(t)\right)
- L\left(t, x(t), u(t)\right)\\
\geq \frac{\partial L}{\partial x}\left(t, \tilde{x}(t), 
\tilde{u}(t)\right)\cdot \left(\tilde{x}(t)- x(t)\right) 
+ \frac{\partial L}{\partial u}\left(t, \tilde{x}(t), 
\tilde{u}(t)\right)\cdot \left(\tilde{u}(t)- u(t)\right)
\end{equation*}
for any control $u$ and its associated trajectory $x$. 
This gives 
\begin{equation}
\label{difrfunctio}
\begin{split}
J[\tilde{x}(\cdot), & \tilde{u}(\cdot)]- J[x(\cdot), u(\cdot)]
= \int^{b}_{a}\left[L\left(t, \tilde{x}(t), \tilde{u}(t)\right)
- L\left(t, x(t), u(t)\right)\right]dt\\ 
& \geq \int^b_a \left[ \frac{\partial L}{\partial x}
\left(t, \tilde{x}(t), \tilde{u}(t)\right)\cdot \left(\tilde{x}(t)- x(t)\right) 
+ \frac{\partial L}{\partial u}\left(t, \tilde{x}(t), 
\tilde{u}(t)\right)\cdot \left(\tilde{u}(t)- u(t)\right)\right]dt\\ 
&= \int^b_a \left[\frac{\partial L}{\partial x}\left(t, \tilde{x}(t), \tilde{u}(t)\right)
\cdot \left(\tilde{x}(t)-x(t) \right) 
- \frac{\partial L}{\partial u}\left(t, \tilde{x}(t), 
\tilde{u}(t)\right)\cdot \left(\tilde{u}(t)-u(t) \right) \right]dt.
\end{split}
\end{equation}
From the adjoint equation \eqref{adj}, we have
\[
\frac{\partial L}{\partial x}(t, \tilde{x}(t), \tilde{u}(t)) 
=  \mathbb{D}^{\psi(\cdot)}_{b^{-}}\lambda(t) 
- \lambda(t)\frac{\partial f}{\partial x}(t, \tilde{x}(t), \tilde{u}(t)).
\]
From the optimality condition \eqref{opt}, we know that 
\[
\frac{\partial L}{\partial u}(t, \tilde{x}(t), \tilde{u}(t))
= - \lambda (t)\frac{\partial f}{\partial u}(t, \tilde{x}(t), \tilde{u}(t)).
\]
It follows from \eqref{difrfunctio} that
\begin{multline}
\label{difrfunctio1}
J[\tilde{x}(\cdot),  \tilde{u}(t)]- J[x(\cdot ), u(\cdot)]
\geq \int^b_a \left(  \mathbb{D}^{\psi(\cdot)}_{b^{-}}\lambda (t) 
- \lambda (t)\frac{\partial f}{\partial x}\left(t, \tilde{x}(t), \tilde{u}(t)\right)\right)
\cdot \left( \tilde{x}(t)-x(t) \right)\\
 - \lambda (t)\frac{\partial 
f}{\partial u}\left(t, \tilde{x}(t), 
\tilde{u}(t)\right) \cdot \left( \tilde{u}(t) -u(t)\right)dt. 
\end{multline} 
Using the integration by parts formula of Lemma~\ref{lemma1}, 
\begin{equation*}
\int^b_a \lambda(t)\cdot ^{C}\mathbb{D}^{\psi(\cdot)}_{a^{+}}\left( 
\tilde{x}(t) -x(t)\right)dt = \left[ \left( \tilde{x}(t) 
-x(t)\right)\cdot \mathbb{I}^{1-\psi (\cdot)}_{b^{-}}\lambda(t) 
\right]^b_a \\
+ \int^b_a \left( \tilde{x}(t) -x(t)\right)
\cdot \mathbb{D}^{\psi(\cdot)}_{b^{-}}\lambda (t)dt,
\end{equation*}
meaning that 
\begin{multline}
\label{integra}
\int^b_a \left( \tilde{x}(t) -x(t)\right)
\cdot \mathbb{D}^{\psi(\cdot)}_{b^{-}}\lambda (t)dt \\
= \int^b_a \lambda(t)\cdot ^{C}\mathbb{D}^{\psi(\cdot)}_{a^{+}}\left( \tilde{x}(t) 
-x(t)\right)dt - \left[ \left( \tilde{x}(t) -x(t)\right)
\cdot \mathbb{I}^{1-\psi (\cdot)}_{b^{-}}\lambda(t) \right]^b_a.
\end{multline}
Substituting \eqref{integra} into \eqref{difrfunctio1}, we get
\begin{multline*}
J\left[ \tilde{x}(\cdot), \tilde{u}(\cdot) \right] 
- J\left[ x(\cdot), u(\cdot)\right]
\geq \int^b_a \lambda(t)\left[ f\left(t, \tilde{x}(t), \tilde{u}(t)\right)\right.
\\ \left.  - f\left(t, x(t), u(t) \right) - \frac{\partial f}{\partial x}\left(t, 
\tilde{x}(t), \tilde{u}(t)\right)\cdot \left( \tilde{x}(t)- x(t)\right) 
- \frac{\partial f}{\partial u}\left(t, \tilde{x}(t), \tilde{u}(t)\right)
\cdot \left( \tilde{u}(t)- u(t)\right)\right]dt.
\end{multline*}
Finally, taking into account that $\lambda(t)\geq 0$ and $f$ is concave 
in both $x$ and $u$, we conclude that 
$J\left[ \tilde{x}(\cdot), \tilde{u}(\cdot) \right] 
- J\left[ x(\cdot), u(\cdot)\right]\geq 0$.
\end{proof}

\begin{Example}
\label{propost}
The extremal $(\tilde{x}, \tilde{u}, \lambda)$ given in Example~\ref{propos} is 
a global minimizer for problem \eqref{pexmple}. This is easily checked 
from Theorem~\ref{theosuff} since the Hamiltonian defined in \eqref{hami} 
is a concave function with respect to both variables $x$ and $u$ and, 
furthermore, $\lambda(t) \equiv 0$. In Figure~\ref{solution_plot}, 
we give the plots of the optimal solution to problem \eqref{pexmple}. 
\begin{figure}[!htb]
\centering 
\includegraphics[scale=0.9]{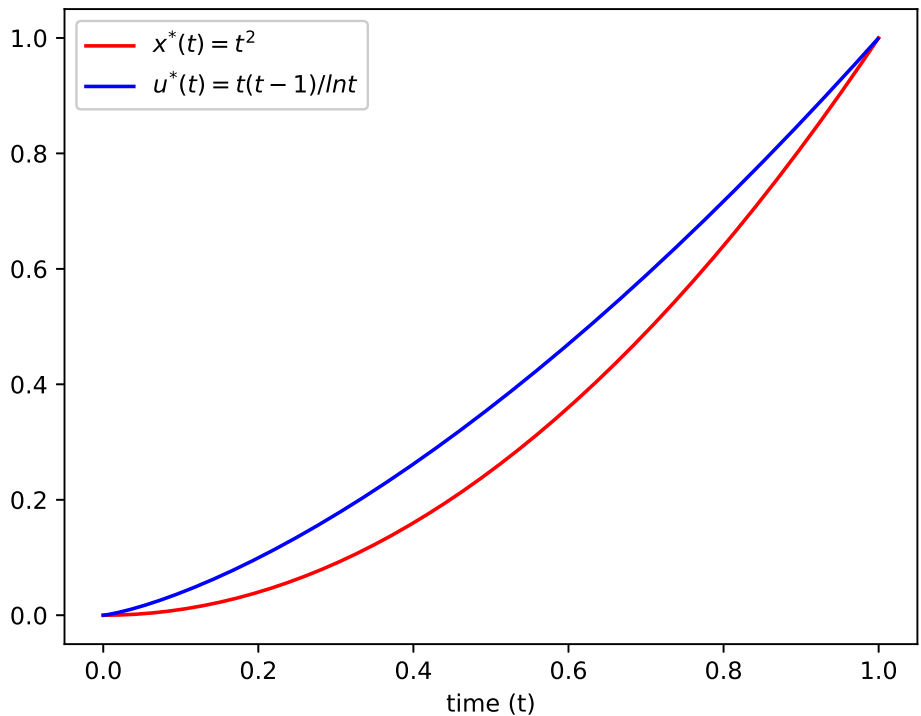}
\caption{The optimal control $u^{*}$ and corresponding optimal state variable $x^{*}$,
solution of problem \eqref{pexmple}.}
\label{solution_plot}
\end{figure}
\end{Example}


\section{Conclusion}
\label{sec:conc}

In this paper we investigated fractional optimal control problems depending 
on distributed-order fractional operators.  We have proved a necessary
optimality condition of Pontryagin's type and a Mangasarian-type sufficient 
optimality condition. The new results were illustrated with an example.
As future work, it would be interesting to develop proper numerical approaches
to solve problems of optimal control with distributed-order fractional
derivatives. In this direction, the approaches found in
\cite{MR3654793} and \cite{zbMATH06915310}
can be easily adapted.


\authorcontributions{The authors equally contributed 
to this paper, read and approved the final manuscript:
Formal analysis, Fa\"{\i}\c{c}al Nda\"{\i}rou and Delfim F. M. Torres; 
Investigation, Fa\"{\i}\c{c}al Nda\"{\i}rou and Delfim F. M. Torres; 
Writing -- original draft, Fa\"{\i}\c{c}al Nda\"{\i}rou and Delfim F. M. Torres; 
Writing -- review \& editing, Fa\"{\i}\c{c}al Nda\"{\i}rou and Delfim F. M. Torres.}

\funding{This research was funded by the Portuguese 
Foundation for Science and Technology (FCT),
grant number UIDB/04106/2020 (CIDMA).
Nda\"{\i}rou was also supported by FCT 
through the PhD fellowship PD/BD/150273/2019.} 

\acknowledgments{The authors are grateful to two anonymous reviewers 
for several comments and suggestions that have helped them to improve 
the manuscript.}

\conflictsofinterest{The authors declare no conflict of interest.} 


\reftitle{References}


\end{document}